\documentstyle{amsppt}
\NoRunningHeads

\topmatter

\title 
On some basic properties of Leibniz algebras
\endtitle

\author   V.Gorbatsevich
\endauthor

\abstract
This paper gives an overview of some basic properties of Leibniz algebras. Some of the results were known earlier, but in the article they are accompanied by new simple proofs. Some of the results are new. The article can be viewed as a digest or a mini-manual for the basic theory of Leibniz algebras 
\endabstract 

\endtopmatter

\document

The main topic of this article is a liezation of Leibniz algebras (the procedure of transition to Lie algebras), which allows to prove many of the properties of Leibniz algebras with aid of known results about Lie algebras. We prove some results recently obtained by various authors. Most of the results below are direct analogs of known results for Lie algebras and almost all of them are proven in many recent articles on the theory of Leibniz algebras (we do not give detailed references here, as the results are very scattered and not always published in the print media, some of them - only in the electronic form). New in this paper is a unified approach and the allocation of some fundamental statements about Leibniz algebras. The article can be considered as a sort of digest or a mini-tutorial on the basics of the theory of Leibniz algebras.

We will consider finite dimensional algebras $L$ over a field $k$ of characteristic 0 (in fact it is only important that this characteristic is not equal to 2). Linear mapping $D: L \to L$ is called a derivation of $L$, if $D (x \cdot y) = Dx \cdot y + x \cdot D (y)$ for all $ x, y \in L $ (where "$ \cdot $" denotes the multiplication in $L$). The set of all derivations of $L$ is denoted by $Der (L)$. Due to the operation of commutation of linear operators $Der (L)$ is a Lie algebra. Naturally can be defined the automorphism group $ Aut (L)$ of the algebra $L$. If the field $k$  is $\bold R $ or $ \bold C $, then the automorphism group is a Lie group and the Lie algebra $Der(L)$ be its Lie algebra. One can consider $Aut (L)$ as an algebraic group (or as a group of $k$-points of an algebraic group defined over a field $k$).

An algebra $L$ over $k$ is called a left Leibniz algebra if for every $x \in L $ the corresponding operator $l_x$ of left multiplication is a derivation of $L$, i.e. $ L_x \in Der (L)$ (Leibniz axiom). Similarly we define the right Leibniz algebra - here operators $r_x$ of right multiplication  have to be derivations of $L$. Leibniz algebra is a generalization of the notion of Lie algebras, which are both left and right Leibniz algebra and their multiplications are skew. So we can say that the Leibniz algebra is a "nonanticommutative" analog of Lie algebra. The term "nonanticommutative" looks very cumbersome (and contains a double negative), sometimes it is used the term "non-commutative" analogue of Lie algebras, although the term is inaccurate - it implicitly assumes that the main operation in the Lie algebras is commutative.

Leibniz algebra at first introduced and investigated in the papers of A.M.Bloh [1] [2] (he called this algebras as $D-$ algebras, directly pointing to their connection with derivations). They have not received a worthy continuation immediately after writing. Then Leibniz algebras were independently rediscovered by Jean-L.Loday (12.01.1946 - 06.06.2012, drowned in the sea riding on a yacht) [3], who called them the Leibniz algebras, since it was Leibniz who proved "Leibniz rule" for differentiation of functions (although there is no information about any "Leibniz algebra" in the works of Leibniz). In the works of Loday (partially - with co-authors), the theory of Leibniz algebras got developed. Some authors call these algebra as Loday's algebras, although the Loday (who sometimes wrote under a pseudonym "Guillaume William Zinbiel", here Zinbiel - is a surname Leibniz, read in reverse order) in [4] pointed that it is a wrong  name. An important role in drawing attention to the Leibniz algebra was played by paper [5] (which was followed by many other works).

A large number of properties (including the fundamental one's) of Leibniz algebras scattered in articles and notes of various authors. While some authors consider the left Leibniz algebra, others - the right one's, making it difficult to read one of these articles. In this paper a listing of some basic definitions and statements of Leibniz algebras is given. However, most of them - with proofs, built on the notion of liezation - a transition to Lie algebras. It turns out that the use of the known properties of Lie algebras dramatically simplify many published proofs in the theory of Leibniz algebras.

The concepts of left and right Leibniz algebras are parallel. For example, we may pass from the right to the left Leibniz algebra by considering a new multiplication $ x \circ y = y \cdot x $. It will be convenient for me to consider left Leibniz algebras, since they have more visible relationship with  differentiation of products (in which the differential operator is written to the right of a differentiable object). However, in many studies on Leibniz algebra the authors prefer to consider right Leibniz algebra. In this case, the main identity is usually written in the form $ x \cdot (y \cdot z) = (x \cdot y) \cdot z - (x \cdot z) \cdot y $. This identity is of course equivalent to the identity, which expresses the property of differentiation for right multiplication by $z$, but to me it seems less natural than based on the concept of differential algebra. In addition, in the study of Leibniz algebras the multiplication is often written (by analogy with the Lie algebras) in bracket form $ [x, y] $. This article will be used in the future such a symbol of multiplication in Leibniz algebras.

For left Leibniz algebras the set of all left multiplication $l_x$ form a Lie algebra. We denote this algebra as  $ad^l(L)$ (similarly for right Leibniz algebras).

Let us consider some simple examples of Leibniz algebras - Leibniz algebras of dimension 1 and 2.

Let $dim_k (L) = 1$ and $a$ - a non-zero element in $ L $. If $[a, a] = 0$, then $L$ is an abelian Lie algebra. If $[a, a] \ne 0$, then $[a, a] = \alpha a$ for some nonzero $ \alpha \in k $. But then the Leibniz identities (both right and left) gives us a contradiction. Therefore, there is only one-dimensional Leibniz algebra - and it is a Lie algebra with the trivial (zero) multiplication.

Now let $dim_k (L) = 2$. Then it is easy to verify that there are only four (up to isomorphism, of course) Leibniz algebra:

1. Two Lie algebras: $a_2$ (abelian) and solvable $ r_2$ (given in a suitable basis $a, b$ by the relations $[a, b] = - [b, a] = b$).

2. Two (left) Leibniz algebra, which are not Lie algebras. In a suitable basis $a, b$ they are given by such relations (we specify only the cases of non-zero products):

(i) $ [b, b] = a $

(ii) $[b, a] = a, [b, b] = a$

Note that the algebra (i) is the left and right Leibniz algebra (obviously), and (ii) only the left one - you can check it by direct calculation of identities. But is easier to use a consequence of Leibniz algebras given below. In fact, for our algebra we have $[[b, b], b] = [a, b] = 0$ (as it should be for the left Leibniz algebra), but $[b, [b, b]] = [ b, a] = a$, but for the right Leibniz algebra it should be $0$. So this two-dimensional left Leibniz algebra is not a right Leibniz one.

The identity of left Leibniz algebra follow some useful consequences of this identity. Here are some of them:

$ r_ {[a, b]} = r_b \circ r_a + l_a \circ r_b $ (for the identity $ [c, [a, b]] $, taking $c$ as an argument)

$ r_ {[ab]} = l_a \circ r_b - r_b \circ l_a $ (from the identity for $ [a, [c, b]] $, taking $c$ as an argument)

These identities imply that $ r_b \circ r_a = - l_b \circ r_a $, i.e. for any $a, b, c \in L$, we have $ [a, b], c]] = - [[b, a], c] $. In particular, for any $a,b$, we have $[[a, a], b] = 0 $. The right Leibniz algebras have the identity of the same type - just rearrange the relevant factors. In particular, for the right Leibniz algebra we have $ [a, [b, b]] = 0$.

Here is a direct proof of the same identities. Take the main identity and the second one obtained from the first by permutation of the arguments:

$ [a, [b, c]] = [[a, b], c] + [b, [a, c]] $

$ [b, [a, c]] = [[b, a], c] + [a, [b, c]] $ or $ [a, [b, c]] = [b, [a, c ]] - [[b, a], c] $

A comparison of these identities gives us $ [[a, b], c] = - [[b, a], c] $.

In the other words, although the operation $ [-, -] $ in the left Leibniz algebra is not skew one, but becomes so after right multiplication by any element. A similar can be said about the right Leibniz algebra - there will need to use multiplication by an element on the left.

Let $ L $ be any left Leibniz algebra. Consider a subset consisting of elements of the form $ [a, a] $ for all possible $ a \in L $. Let us consider the linear subspace generated by this subset (we denote it by $ Ker (L) $; reason for such notation and other possible its name age written below) - it is an ideal in $ L $, left and right one's, abelian and lefty-center (the multiplication $ [Ker (L ), L] $  is 0.) The fact that it is the right ideal follows from such identity: $ [a, [b, b]] = [a + [b, b], a + [b, b]]] - [a, a] $. But the ideal $ Ker (L) $ is not always right-central. This ideal (using other notation) was firstly introduced in [2].

Here is an another proof of the fact that $ Ker (L)$ is a right ideal: at first we show that $ Ker (L) $ is spanned by elements of the form $ [a, b] + [b, a] $ (this follows from the identity $ [a,b] + [b,a] = [a + b, a + b] - [a, a] - [b, b] $), and then we use the identity $ [a, [b, b]] = [c, b] + [b, c] $, where $ c = [a, b] $.

The ideal $ Ker (L) $ always is different from the $ L $, as it abelian and abelian Leibniz algebra is a Lie algebra and therefore for it $Ker (L) =  \{ 0 \} $.

Factor algebra $ L / Ker (L) $ is a non-zero algebra, where $ Ker (L) $ is the smallest ideal in $ L $, the quotient algebra by which is a Lie algebra (both for left and for right Leibniz algebras). Therefore $ Ker (L)$ can be viewed as a "non-Lie core" and perhaps could be called the "liezator" for $ L $. The  ideal $ Ker (L) $ can be described in another way - as the linear span of all elements of the form $ [x, y] + [y, x] $. Factor algebra $ L / Ker (L) $ is called as the liezation of $ L $ and denotes by $ L ^ \star $. There is a natural action of the Lie algebra $ L^\star $ on a vector space $ Ker (L) $ (multiplication in which is trivial). 

If the Leibniz algebra $ L $ is commutative (i.e. $ [x, y] = [y, x] $), then the subset $ [x, y] + [y, x] $ coincides with the commutant $ [L, L ] = \{ [x, y] | x, y \in L \} $. But then liezation of the algebra $ L $ is just an abelian Lie algebra $ L / [L, L] $. Commutative algebras are nilpotent Leibniz algebras (about which see the next section), their class of nilpotency equals to 2. This kind of Leibniz algebra can be described in some detail.

Now let us consider the left center $Z^l (L) = \{ x \in L | [x, L] = 0 \} $. Similarly we introduce the concept of the right center $Z^k (L)$ (both this centers can be considered for the left and right Leibniz algebras). For the left Leibniz algebra $ Z^l (L) $ is an ideal, moreover - two-sided ideal (since $ [[x, y], L] = - [[y, x], L] $). The right center $Z^k ( L) $ is an  subalgebra (because $[x,[u,v]] = [[x,u],v]+[u,[x,v]] = 0$ for all $x \in L, u,v \in Z^r(L)$, and in general, the left and right centers are different; they even may have different dimensions. By the above, we have $ Ker (L) \subset Z^l (L) $. Therefore $ L / Z^ l (L) $ is the Lie algebra, which is isomorphic to the Lie algebra $ad^l(L)$ mentioned above.

For Leibniz algebras may be naturally introduced many of the concepts studied in the theory of Lie algebras. For example, consider series (derived series) of commutants $ D^ n (L) $: $ D^1 (L) = [L, L], D^{k +1} (L) = [D^k (L), D^k (L)] $. It is important to note that the members of this series are two-sided ideals (both for  the right and for the left Leibniz algebras) - it is an easy consequence of the Leibniz identity. Leibniz algebra is called solvable if the derived series comes to $ \{0 \} $ at some finite step. It is easy to check that the sum of solvable ideals in Leibniz algebra is also a solvable ideal. Therefore, in $L$ there exists a largest solvable ideal $ R $ (containing all other soluble ideals). Naturally, it is called the radical of the Leibniz algebra $ L $. Since $ Ker (L) $ is an ideal with trivial multiplication, then it is contained in the radical of $ R $ of Leibniz algebra $ L $.

Next, we consider the decreasing central series for Leibniz algebra $L$: $ C^ n (L)$: $ C^ 1 (L) = [L, L]$, $C^ {k +1} (L) = [L, C^ k (L)] $. Despite a certain lack of symmetry of the definition (multiplication by $ L $ only from the left) members of this series are two-sided ideals, since by the Leibniz identity (any - right or left), we have $ [C^ k (L), L] = [L, C^ k (L)] $. Therefore, $ C^k (L) $ consists of linear combinations of the products of $k$ elements with arbitrarily spaced brackets. In particular, it implies such inclusion $[C^p(L),C^q(L)] \subset C^{p + q}(L)$.

Leibniz algebra is called nilpotent if its central series reaches zero in a finite number of steps. As it is follows from definition the centers (left and right one's) are nonzero.

\proclaim {Proposition 1} Any Leibniz algebra $L$ has a maximal nilpotent ideal (containing all the nilpotent ideals of $L$)
\endproclaim

\demo {Proof}
Let $I$ be some nilpotent ideal in $ L $. Then the sum $ I + Ker (L) $ is a nilpotent Ideal too - this follows immediately from the left centrality of $Ker (L) $. Therefore, in discussing the question of the existence of maximal nilpotent ideal is sufficient to consider only nilpotent ideals containing $ Ker (L) $.

Consider the epimorphism of liezation $ \Cal L: L \to L / Ker (L) = L^\star $. For the Lie algebra $ L / Ker (L)$ the existence of a maximal nilpotent ideal is known. Its inverse image under the mapping $L \to \Cal L $ also is nilpotent (by the centrality of $ Ker (L) $) and it is maximal because of what is said above about the maximum nilpotent ideals in $ L $.
\enddemo

By nilradical in a Leibniz algebra $L$ is called a maximal nilpotent ideal in $L$ (which exists by Proposition 1). By its definition, the nilradical is a characteristic ideal, i.e. It remains invariant under all automorphisms of the Leibniz algebra $L$. It is obviously contained in the radical of the Leibniz algebra and it equals to the nilradical of the solvable radical of $L$. The nilradical contains right center, as well as the ideal $Ker(L) $.

We also consider the notion of normalizer. The left normalizer $ N^ l_L (U) $ of subset $ U \subset L $ in Leibniz algebra $ L $  is the set of such $x  \in L $, that $ [a, U] \subset U $. For the right normalizer $ N^r_L (U) $ it requires that $ [U, a] \subset U $. For the left Leibniz algebra a left normalizer of some subalgebra $M$ is a subalgebra, but the right normalizer of $M$ is not necessarily a subalgebra (but see Lemma below).

\proclaim {Proposition 2} (Engel's theorem for Leibniz algebras) If all operators $ l_x $ of left multiplication for the left Leibniz algebra $ L $ are nilpotent, then the algebra $ L $ is nilpotent In particular, for left multiplications there is a common eigenvector with zero eigenvalue

In some basis matrices of all $ l_x $ are upper-triangular nilpotent matrices.
\endproclaim

\demo {Proof} For the Lie algebra $ L / Ker (L) $ the Engel's theorem is well known. Therefore, the assertion of Proposition 2 holds for it. But $ Ker (L) $ is the central ideal - it action by left multiplication on $ L $ is trivial. Therefore Leibniz algebra $ L $ can be considered as a central extension of a nilpotent Lie algebra - so it is also nilpotent.

Reducing of matrix of left multiplications to upper-triangular nilpotent form is obvious now (using induction).
\enddemo

There is a stronger version [6] of this proposition: if for left Leibniz algebra all left
multiplication are nilpotent, then all right multiplications also are nilpotent and in a suitable bases linear operators of all left and all right multiplications have upper-triangular nilpotent matrices. Moreover, both left and right multiplications have a common eigenvector with eigenvalue equals to 0.

The proof of this version is based on the identity $ (r_x)^ n = (-1)^n r_x (l_x)^{n-1} $, which is proved by induction. For $ n = 2 $, it follows from the identity mentioned above: $ [[x, a], a] = - [[a, x], a] $. Then take $ (r_x)^3 (a) = [[[a, x], x], x] $ and transform it similarly and then we may use an induction.

\proclaim {Corollary 1} The normalizers (both - left and right) of some subalgebra $ M $ in a nilpotent Leibniz algebra $ L $ are not equal to the subalgebra $ M $ (they both  strictly contain $M$)
\endproclaim

\demo {Proof} At first we consider the left normalizer. The reasoning here is exactly the same as for Lie algebras - we consider the action induced by left multiplication on the factor space $ L / M $. The set of linear operators of left multiplications forms a Lie algebra $ad^l(L)$, which is nilpotent in our case. So by Engel's theorem for Lie algebras this action has a common eigenvector with 0 as an eigenvalue. The vector in $L$ which corresponding to this eigenvector is contained in $ N_L (M) $ and is out of $ M $. For the right normalizer we use the stronger version of Engel theorem mentioned above.
\enddemo

\proclaim {Corollary 2} A subspace  $ V \subset L $ generates a Leibniz algebra if and only if $ V + [L, L] = L $
\endproclaim

\demo {Proof}
Let $ V + [L, L] = L $. As $ M $ we denote the subalgebra of $ L $, generated by the subspace $ V $. We want to prove that $ M = L $. The proof uses inductions on $codim(M)$ and $dim(L)$. For $codim (M)=0 $ the statement is obvious, because in this case $ M = L $. For $ dim (L) = 0 $ this statement is evident too.

Let $ codim (M)> 0 $. Consider the normalizer $ M^ \prime = N_L (M) $ - it is strictly greater than $ M $. Therefore $ codim (M^ \prime) < codim (M)$ and so by the induction we have $ M^ \prime = L $. This means that $ M $ is the left ideal in $ L $, and it is supplemental to the $ [L, L] $. Now we are going to prove that $M$ may be supposed to be a right ideal too. We use one easy lemma.

\proclaim{Lemma} Let $I$ be a left ideal in a left Leibniz algebra $L$. Then the right normalizer $N^r_L(I)$ is a subalgebra in  $L$.
\endproclaim

\demo{Proof of Lemma} Let us consider two elements $n, n^\prime \in N^r_L(M)$. We have $[M,[n,n^\prime]]$ = $[[M,n], n^\prime] +[n, [M,n^\prime]]$. But $[M,n]$ and $[M,n^\prime]$ both contain in $M$ (due to the definition of the normalizer) and $[n,M] \subset M$ because $M$ is the left ideal in $L$. So $[M, [n,n^\prime]] \subset M$,
$[n,n^\prime] \in N^r_L(M)$, i.e. $M$ is the subalgebra in $L$. 
\enddemo

In our case $M$ is the left ideal in $L$. Due to Lemma $N^r_L(M)$ is a subalgebra in $L$ (and it strictly contains $M$ - see Corollary 1). We continue to use an induction on $codim_L(M)$. We have $codimL(N^r_L(M)) < codim_L (M)$. Therefore  by induction we have  $N^r_L(M) = L$, i.e. $M$ is the right ideal in $L$.

Therefore we may suppose now that $M$ is a two-sided ideal in $L$.

We use now a induction on $dim (L)$. Let $ L^\star = L / Ker(L) $ and consider the natural epimorphism $ \Cal L: L \to L^\star $. It is clear that $ \Cal L ([L, L]) = [L^\star, L^\star] $. Let $ V^\star = \Cal L (V) $. It is clear that $ V^\star + [L^\star, L^\star] = L^\star $. Since the dimension of the Lie algebra $ L^\star $ less than the dimension of the Lie algebra $ L $, then we can hold the induction step. Therefore, we conclude that the subspace $ V^\star $ generates a Lie algebra $ L^\star $. But then we have $ M + Ker(L) = L $. 

For the commutant of $ L $ we have $ [L, L] = [M + Ker(L), M + Ker(L)] = [M, M] + [M, Ker(L)] + [Ker(L), M] + [Ker(L), Ker(L)] $. However, the subspaces $ [Ker(L), Ker(L)] $ and $ [Ker(L), M] $ are zero by left centrality of $ Z $. Also we have $ [M, Ker(L)] \subset M $, since $ M $ also is the right ideal in $ L $. So $ [L, L] = [M, M] + M \subset M$, and therefore $ M $ contains $ [L, L] $. Since $ M + [L, L] = L $, we see that $ M = L $.

Conversely, if the subspace of $ V $ generates a Leibniz algebra $ L $, then its image under the natural epimorphism $ L \to L / [L, L] $ generates abelian Leibniz algebra $ L / [L, L] $. Since the multiplication in $ L / [L, L] $ is trivial, then the subspace of this algebra can generate it only if this subspace coincides with the whole algebra. But this is equivalent to $ V + [L, L] = L $
\enddemo

It is interesting to note that not all properties of nilpotent Lie algebras, even a simple and well-known one's, are hold for the case of Leibniz algebras. For example, there is a simple statement for nilpotent Lie algebras of dimension 2 or more: "codimension of commutant is $ \ge 2 $". For Leibniz algebras it is not true (though not only for nilpotent, but for all solvable Leibniz algebras we have $ codim_L [L, L]> 0 $). For example, two-dimensional Leibniz algebra $ <a,b> $, for which $ [a, a] = b $ (see above) is nilpotent, but its commutant has codimension 1. This is due to the fact that its liezation is one dimensional. But for one-dimensional Lie algebras above mentioned statement is incorrect. 
  
We obtain the useful corollary.

\proclaim {Corollary 3} If Leibniz algebra $ L $ is nilpotent and $ codim_L ([L, L]) = 1 $, the algebra $ L $ generated by one element.
\endproclaim

So for $ codim_L ([L, L]) = 1 $ a nilpotent algebra is a kind of "cyclic". The study of such nilpotent algebras is the specifics of the theory of Leibniz algebras; Lie algebras has no analogue. Such cyclic $L$ may be explicitly described.

\proclaim {Corollary 4} The minimal number of generators of Leibniz algebra $L$ equals to $ dimL / [L, L] $
\endproclaim

We recall now that a complete flag of subspaces of a finite dimensional vector space $ V $ is a family of nested subspaces $ V_0 = \{ 0 \} \subset V_1 \subset \dots \subset V_n = V $ such that $ dim V_k = k $.

\proclaim {Proposition 3} (weak version of Lie Theorem for solvable Leibniz algebras) A left Leibniz algebra $ L $ over $ \bold C $ has a complete flag of subspaces which is invariant under left multiplication.

In other words, all linear operators $ l_x $ of left multiplications can be simultaneously reduced to triangular form.
\endproclaim

\demo {Proof} For Lie algebras this theorem in well known (in one of the variants  of its formulation). Now let $ L $ be an arbitrary solvable Leibniz algebra. Since liezator $ Ker (L) $ is central, then the solvability of the Leibniz algebra $ L $ is equivalent to the solvability of the Lie algebra $ L / Ker (L) $. By the classical Lie theorem $ L / Ker (L) $ has a complete flag, which is invariant under the multiplication (both left and right).

Let $ x \in L $ be some element. Its left action on $ Ker (L) $ is trivial and therefore in $ Ker (L) $ there is a complete flag, which is invariant under left multiplication (for the right multiplication on the left Leibniz algebra such reasoning does not pass). Adding to it the full invariant flag in the factor space until the full flag in $L$, we obtain a complete flag in $ L $, which is invariant under left multiplication.
\enddemo

\proclaim {Proposition 4} Let $ R $ - radical in the Leibniz algebra $ L $, and $ N $ - nilradical of $ L $. Then $ [L, R] \subset N $
\endproclaim

\demo {Proof}
For Lie algebras $ L $ that statement is true. Now let $ L $ be an arbitrary Leibniz algebra. Then $ Ker (L) \subset N \subset R $. Let $ L^ \star = L / Ker (L) $ (the liezation of algebra $ L $), $ R^ \star = R / Ker (L), N^ \star = N / Ker (L) $. From the definitions of the radical and nilradical it follows that $ R^ \star $ $ and N^ \star $ are radical and nilradical of Lie algebra $ L ^ \star $ respectively (we remember that $ Ker (L) \subset Z^l (L ) $). Since $ L^ \star $ is the Lie algebra, then $ [L^ \star, R^ \star] \subset N^ \star $. But then, by the fact that $ Ker (L) \subset N \subset R $, we get $ [L, R] \subset N $
\enddemo

I note that the statements contained in Propositions 3 and 4 are not met me in the papers on the theory of Leibniz algebras (although they are well known for Lie algebras, especially Lie Theorem).

\proclaim {Corollary 5} $ [R, R] \subset N $. In particular, $ [R, R] $ nilpotent.
\endproclaim

\proclaim {Corollary 6} Leibniz algebra $ L $ is solvable if and only if $ [L, L] $ is nilpotent.
\endproclaim

\demo {Proof} In one direction the assertion of Corollary 6 is proved in Corollary 5. Let us consider the converse.

Let $ [L, L] $ be nilpotent. Consider the radical $ R $  of algebra $ L $. The Lie algebra $ L / R $ is the semi-simple Lie algebra (as it contains no solvable ideals). But semisimple Lie algebra coincides with its commutator. Therefore, in the case where $ [L, L] $ is nilpotent, the this semisimple algebra should be trivial. But this means that $ L $ is solvable.
\enddemo

\proclaim {Proposition 5} Multiplications (right and left one's) in the Leibniz algebra are degenerate linear operators (their determinants equal to 0).
\endproclaim

\demo {Proof} For Lie algebras this degeneration is obvious, since $ [x, x] = 0 $
for any $ x \in L $, i.e. the vector $ x $ is the root vector for the zero eigenvalue.

     In any Leibniz algebra $ Ker (L) $ is the two-sided ideal different from $ L $ (see above) and $ L / Ker (L) $ is a Lie algebra of positive dimension. Let $ x $ be some non-zero element in $ L $. For $ l_x $ and $ r_x $ the subspace $ Ker (L) $ invariant. On the factor space which is  a Lie algebra, the linear operators induced by these multiplications, are degenerated (their determinants equal to 0). But then these multiplications are degenerate in the whole space $ L $.
\enddemo

There is an interesting question about the existence of faithful representations of Leibniz algebras. Linear representation (sometimes referred as module) of Leibniz algebra is a vector space $ V $, for which we have two actions (left and right) of the Leibniz algebra $ L $ $ [-, -]: L \times M \to V $ and $ [- , -]: M \times L \to V $, such that the identity

$ [x, [y, z]] = [[x, y], z] + [y, [x, z]] $

is true if one (any) of the variables is in $ M $, and two others - in $ L $. In some other formulation (and without using bracketed notation for the action of the algebra $ L $ on $ V $) we have such conditions conditions:

$ a(bm) = [ab]m + b(am) $

$ a(mb) = (am)b + m[ab] $

$ m[ab] = (ma)b + a(mb) $

Note that the concept of representations of Lie algebras and Leibniz algebras are different. Therefore, such an important theorem in the theory of Lie algebras, as the Ado theorem on the existence of faithful representation in the case of Leibniz algebras was proved much easier and gives a stronger result. It is because the kernel of the Leibniz algebra representation is the intersection of kernels (in general, different one's) of right and left actions, in contrast to representations of Lie algebras, where these kernels are the same. Therefore, an faithful representation of the Leibniz can be obtained easier than faithful representation of the Lie algebra.

\proclaim {Proposition 6 (see [7])} Any Leibniz algebra has a faithful representation of dimension $ \le dim (L) +1 $.
\endproclaim

The structure of Leibniz algebras may be described by such result. 

\proclaim {Proposition 7 (Levi theorem for Leibniz algebras, see [8])} For a Leibniz algebra $L$ there exists a subalgebra $ S $ (which is a semisimple Lie algebra), which gives the decomposition $ L = S + R $, where $ R $ - radical.
\endproclaim

\demo {Proof} Consider the Lie algebra $ L^ \star = L / Ker (L) $. By the classical Levi theorem there is a semisimple algebra $ S^\star \subset L^ \star $, which gives the semidirect sum decomposition $ L^\star = S^\star + R^\star $. Here $ R^\star $ is the radical of $ L^ \star $. Let  $ R $ be the radical of Leibniz algebra $ L $ and $ \Cal L: L \to L ^ \star $ be the natural epimorphism. Then it is easy to understand, that $ \Cal (R) = R^\star $. Let $ F = \Cal L^{-1} (S^ \star) $. Then $ F $ is the  Leibniz subalgebra in $ L $ which contains $ Ker (L) $, and the quotient algebra $ F / Ker (L) $ is isomorphic to a semisimple Lie algebra $ S^\star $.

      For the Lie algebra $ S^ \star $ the abelian Leibniz algebra $ Ker (L) $ is $ S^ \star-$ module. For semisimple Lie algebra $ S^ \star $, it follows that the subalgebra $ F $, considered as an extension of $ S^ \star $ using $ Ker (L) $, is split (by Whitehead's Lemma for semisimple Lie algebras). But this means that in $F$ there exists a subalgebra (semisimple Lie algebra) which is complementary to the abelian ideal $ Ker (L) $. Clearly, the $ S + R = L $ and we get the required decomposition.
\enddemo

By examples (see [8]) the non-uniqueness of the subalgebra $ S $ it can be shown  (minimum dimension of Leibniz algebra there equals to 6). In the case of Lie algebras the semi-simple Levy factor is unique up to conjugation.

\proclaim {Corollary 7} Let $L$ be a Leibniz algebra such that its left center lies in the right center (i.e., that $ [L, Z^l (L)] = \{ 0 \} $). Then if the Lie algebra $ L / Z^l (L) $ is a semisimple one, that $ L $ is the reductive Lie algebra (so its radical is central)
\endproclaim

\demo {Proof} Let us consider the Levi decomposition $ L = S + R $. Since, by hypothesis, the algebra $ L / Z^ r (L) $ is semisimple, then $ R = Z^l (L) $. The left action of the subalgebra $ S $ on $ Z^l (L) $ is always trivial. But since $ [L, Z^ l (L)] = \{ 0 \} $, the right action of $ S $ on $ Z^ l (L) $ is also trivial. Therefore, $ S $ commutes with the radical and therefore Levi decomposition of $ L $ is a direct one. The radical of $ R = Z^l (L) $ is abelian, so $ L $ is the reductive Lie algebra.
\enddemo

Statement of Corollary 7 is proved in [5] for the special case when $ S $ is a  classical simple Lie algebra.

The author is grateful to B.Omirov  for helpful discussions about the theory of Leibniz algebras.

\enddocument